\def\fs{\footnotesize}
\def\bm{\begin{em}}
\def\en{\end{em}}
\def\be{\begin{equation}}
\def\ee{\end{equation}}

\def\cd{\cdot}

\def\ra{\rightarrow}
\def\C{{\cal C}}
\def\bk{I \kern-.25em k} 
\def\C{{\cal C}}
\def\nd{\hbox{\rm End}\,}
\def\hom{\hbox{\rm Hom}\,}
\def\com{\hbox{\rm Com}\,}
\def\id{\hbox{\rm id}}
\def\vec{{\cal V}{\rm ect}_{\bk}}
\def\de{\delta}
\def\ep{\varepsilon}
\def\al{\alpha}
\def\bt{\beta}
\def\ga{\gamma}

\def\O{\mit\Omega}
\def\sw{\rm SW}

\def\co{\triangle}
\def\cl{\bigtriangleup_L}
\def\ccr{\bigtriangleup_R}
\def\clb{\bar\bigtriangleup_L}
\def\crb{\bar\bigtriangleup_R}
\def\bu{\bullet}
\def\da{\dagger}

\def\ot{\otimes}
\def\c{\circ}

\documentclass[12pt]{article}
\usepackage{latexsym}
\begin{document}
\title{SOME TOPICS IN COALGEBRA CALCULUS}
\date{Czech. J. Phys. {\bf50} no1 (2000), 23-28}
\author{
Andrzej Z. BOROWIEC\thanks{
Email: borow@ift.uni.wroc.pl .}\\
{\fs Institute of Theoretical Physics, University of Wroc{\l}aw}\\
{\fs pl. Maksa Borna 9, PL-50-204  WROC{\L}AW, POLAND}\\
Guillermo A. V\'AZQUEZ--COUTI\~NO\thanks{
E-mail: gavc@xanum.uam.mx .}\\
{\fs Universidad Aut\'onoma Metropolitana-Iztapalapa}\\
{\fs A.P. 55-534, M\'EXICO D.F., C.P. 09340}\\
}
\maketitle

\begin{abstract}
 We study a coderivation from a cobimodule into a coalgebra. Vector cofields
are defined by the action of a codual bicomodule on  a coalgebra. This
action is induced by a codifferential. A construction of a codual object
in the category of bicomodules over a given coalgebra has been proposed.
Various notions of duality have been analysed in this context.
\end{abstract}

\section{Introduction}

It is well known, that the notion of coalgebra and comodule is dual
to the concept of algebra and module: the dual vector space to a linear
finite dimensional algebra inherits a canonical coalgebra structure
(see e.g. \cite{Abe,SW}). Similarly, dualizing the concept of derivation one
can obtain the notion of coderivation, etc.

A general  vector field  formalism for a first order differential
calculus over a given algebra $A$ has been recently proposed
in \cite {AB1,AB2}. Our aim in the present note is to introduce a dual
concept of vector {\it cofields} for a codifferential calculus over
a linear coalgebra $\C$.

To clarify the above idea, one should notice that there are several notions
of duality. E.g. the already mentioned duality between algebras and coalgebras,
modules and comodules holds in the category of (finite dimensional)
vector spaces over a fixed field $\bk$ \cite{Abe,SW}. Also one-forms and
vector fields are dual concepts, but this time duality holds in the
category of bimodules over the algebra $A$ \cite {AB1,AB2}. Similarly,
a major step in our program relies upon introducing a notion of duality in the
category of bicomodules over the linear coalgebra $\C$. The present letter has
a preliminary character. More detailed study of this subject will be
presented elsewhere. Some proofs are left to the reader.

As explained above, we shall work mainly in the category $\vec$ of
finite dimensional vector $\bk$-spaces. However, the authors believe
that all results can be extended after adding suitable topological conditions
or/and to the case of {\it rational} comodules (cf. \cite{Abe,SW}).

Throughout this letter we let $\bk$ be a field. All objects considered
here are finite dimensional vector spaces over the field $\bk$ and all
maps are $\bk$-linear maps. Tensor products are also assumed over $\bk$.
Given vector spaces $U$ and $W$, we denote by $\hom (U, W)$ the space
of linear maps from $U$ to $W$. For any space $U$, let us denote by
$U^\bu\doteq\hom (U, \bk)$  the linear dual of $U$.

Due to finite dimensionality the following canonical isomorphisms:
$$
V\ot W\cong W\ot V, \ \ \hom (V, W)\cong V^\bu\ot W, \ \
\hom (V, W\ot U)\cong \hom(V, W)\ot U
$$ will be utilized without special mention.

For the same reason, we found it also appropriate to adopt the covariant
index notation, together with the Einstein summation convention over
repeated up ({\it contravariant}) and down ({\it covariant}) indices,
borrowed from tensor analysis on manifolds and General Relativity.
For example, an element $u\in U$ will be  written as $u=u^i e_i$,
where $\{e_k\}_{k=1}^{{\rm dim}U}$ denotes some basis in $U$.
Thus our results, although formally formulated in a
basis-dependent way, are in fact basis-independent.

Our main references to coalgebras and comodules are \cite{Abe,SW}.
For notational convenience, we fix some coassociative and counital
coalgebra $\C\equiv (\C,\co ,\ep)$ in $\vec$. All comodules will be
considered over the coalgebra $\C$.

\section{Bicomodules and their coduals}

If $(U, \cl)$ is a left comodule with a (finite dimensional) carrier
space $U$ then in an arbitrary basis $\{e_k\}$ for $U$ one can set
\be
\cl (e_k)=L^i_k\ot e_i \ \ \ \ \hbox{(summation convention !)},
\ee
where $\C$-valued matrix elements $L^i_k$ have to satisfy the following
relations (see e.g. \cite{Abe}):
\be
\co(L^i_k)=L^m_k\ot L^i_m ,\ \ \ \ep (L^i_k)=\de^i_k .
\ee
Similarly, a right comodule structure $\ccr$ on $U$ can be encoded by
\be
\ccr (e_k)=e_i\ot R^i_k,
\ee
with
\be
\co(R^i_k)=R^i_m\ot R^m_k ,\ \ \ \ep (R^i_k)=\de^i_k .
\ee
The left comodule structure $\cl$ on $U$ induces by the transpose
coaction a canonical \underline{right} comodule structure $\cl^\bu$
on $U^\bu$. Taking the dual basis $\{e^k\}$ in $U^\bu$ one defines:
\be
\cl^\bu (e^k)\doteq e^m\ot L^k_m= (\id\ot e^k)\c\cl ,
\ee
as equalities in $\hom (U, \C)$. Similarly, $\ccr$ generates
\be
\ccr^\bu (e^k)\doteq R^k_m\ot e^m =(e^k\ot\id)\c\ccr
\ee
a canonical \underline{left} comodule structure on $U^\bu$.\medskip\\
{\bf Remark 2.1}. It should, however, be noticed that the interpretation
of $\cl^\bu$ as a right comultiplication in $U^\bu$ requires the identification
$\hom (U,\C)\cong U^\bu\ot\C$. Similarly, to see $\ccr^\bu$ as a left
comultiplication one has to identify $\hom (U,\C)$ with $\C\ot U^\bu$ instead.\medskip\\

If $(U, \cl , \ccr)$ is a  bicomodule, then the commutation relation
$(\id\ot\ccr)\c\cl=
(\cl\ot\id)\c\ccr$ gives rise to the following identity in $\C\ot\C$
\be
L^i_k\ot R^m_i=L^m_j\ot R^j_k .
\ee
In the above situation $(U^\bu, \cl^\bu , \ccr^\bu)$ becomes
automatically a bicomodule too.\medskip\\
{\bf Example 2.2}. Choosing an arbitrary basis $\{c_\al\}$ in the vector
space $\C$, one can define
\be
\co (c_\al) \doteq \O^{\bt\ga}_\al c_\bt\ot c_\ga \doteq
l^\ga_\al\ot c_\ga \doteq c_\bt\ot r^\bt_\al \ ,
\ee
where $\O^{\bt\ga}_\al\in\bk$ are {\it structure constants} for $\C$;
$l^\ga_\al\doteq \O^{\bt\ga}_\al c_\bt\in\C$ and
$r^\bt_\al\doteq \O^{\bt\ga}_\al c_\ga\in\C$  do define a canonical
bicomodule structure $(\C, \co_l , \co_r)$
with $\co\equiv\co_l\equiv\co_r$ being mappings from $\C$ into
$\C\ot\C$. This means,
that its dual $(\C^\bu, \co^\bu_r , \co^\bu_l)$ is a bicomodule too
(but not a coalgebra !).
Of course, $\C^\bu$ possesses a canonical
(unital, associative) algebra structure with a multiplication table
given by the same structure constants:
\be
m(c^\al\ot c^\bt)\equiv c^\al\cd c^\bt \doteq \O^{\al\bt}_\ga c^\ga \ .
\ee
The counit $\ep\in\C^\bu$ is a unit for this algebraic structure and
\be
\co^\bu_l (\ep) = \co^\bu_r (\ep) = c^\al\ot c_\al=\id_\C .
\ee
Notice, as a side remark, that algebra and comodule structures
on $\C^\bu$ are related by two compatibility conditions:
\be
\co^\bu_r\c m = (\id\ot m)\c (\co^\bu_r \ot\id)\ ,\ \ \
\co^\bu_l\c m = (m\ot\id)\c (\id\ot\co^\bu_l)\ ,
\ee
i.e. $\C^\bu$ carries a double {\it dimodule} structure in the
terminology of \cite{BP}.

We recall a known fact from the theory of modules: if $A$ is a
linear algebra and $M$, $N$ are $A$-bimodules then $\hom (M, N)$
is a quadruple $A$-module. It means that $\hom (M, N)$ has an
$A$-module structure in four different ways: two induced by the
multiplications in $M$ and another two induced by the
multiplications in $N$. Moreover, these different $A$-module
structures are pair-wise commuting.

One needs a comodule analogue of this.
Let us begin with a simple example.\medskip\\
{\bf Example 2.3}. Given two comodules $(U, \cl)$ and $(W, \crb)$.
The vector space $U\ot W\cong W\ot U$ can be endowed with an obvious
bicomodule structure by
\be
\cl^{U\ot W}=\cl\ot\id_W , \ \ \ \ \ \ccr^{U\ot W}=\id_U\ot\crb.
\ee
Let now $(U,\cl,\ccr)$ and $(W,\clb,\crb)$ be two bicomodules.
In the vector space $\hom(W,U)$ one can introduce four different
comodule structures induced correspondingly by $\crb,\ \clb,\ \cl,\ccr$:
\be
\cl^1 (\Phi)=(\Phi\ot\id)\c\crb\in\hom (W,U\ot\C) ,
\ccr^1 (\Phi)=(\id\ot\Phi)\c\clb\in\hom (W,\C\ot U) ;
\ee
\be
\cl^2 (\Phi) = \cl\c\Phi\in\hom (W,\C\ot U) ,
\ccr^2 (\Phi) = \ccr\c\Phi\in\hom (W, U\ot\C) .
\ee
{\bf Remark 2.4}. Again, in order to ensure a proper interpretation for the
corresponding comultiplications, one has to identify
$\hom (W,U\ot\C)\cong\C\ot\hom (W,U)$
and $\hom (W,\C\ot U)\cong\hom (W,U)\ot\C$ in (13).
On the contrary, (14) requires that $\hom (W,\C\ot U)\cong\C\ot\hom (W,U)$
and $\hom (W,U\ot\C)\cong\hom (W,U)\ot\C$ instead.\smallskip\\

According to the Example 2.3, after making the identification
$\hom (W,U)\cong W^\bu\ot~U\cong U\ot W^\bu$ one can also write:
\be
\cl^{1}\doteq\crb^\bu\ot\id_U,\ \ \ \ \ccr^{1}\doteq\id_U\ot\clb^\bu,
\ee
\be
\cl^2\doteq\cl\ot\id_{W^\bu},\ \ \ \ \ccr^2\doteq\id_{W^\bu}\ot\ccr.
\ee
This argue in favour of the following\medskip\\
{\bf Proposition 2.5}. \bm Above four comodule structures are pairwise commuting,
i.e. \\$(\hom (W, U), \cl^1, \ccr^1, \cl^2, \ccr^2)$ is a
quadruple-comodule.\en\\
{\bf Remark 2.6}. In particular, $\nd\C$ is a quadruple-comodule too.
In what follows, we shall use $(\nd\C, \cl^1, \ccr^1)$ as a default
bicomodule structure on $\nd\C$. Thus, for $\xi\in\nd\C$
\be
\cl^1 (\xi) = (\xi\ot\id)\c\co, \ \ \ \ \ccr^1 (\xi)=(\id\ot\xi)\c\co .
\ee\\
\indent If $\Phi\in\hom (W, U)$ is a left comodule map from
$(W, \clb)$ into $(U, \cl)$ then
\be
\cl\c\Phi = (\id\ot\Phi)\c\clb .
\ee
In our notation  it reads $\cl^2(\Phi) = \ccr^1(\Phi)$ in $\hom (W,\C\ot U)$.
Select  the subspace $\com^{(\C,-)}(W, U)\subset\hom (W, U)$ of all left
comodule maps from $W$ into $U$. It turns out that this subspace becomes a
subbicomodule with respect to $(\cl^1, \ccr^2)$. Similarly, the space of right
comodule maps $(\com^{(-,\C)} (W, U), \cl^2, \ccr^1)$ is also a bicomodule.\\

With all these and keeping in mind an anlogous situation in the category
of bimodules over an algebra, described in \cite{AB1,AB2}, we are ready to
introduce a codual object in the category $^\C\!{\cal M}^\C$ of
$\C$-bicomodules. Let $(U,\cl,\ccr)$ be a bicomodule.\medskip\\
{\bf Definition 2.7}.
A bicomodule $(^\da\! U\equiv \com^{(\C,-)}(\C,U),\cl^1,\ccr^2)$ is
called a left $\C$-codual of $U$. In a similar manner, one defines
$(U^\da\equiv\com^{(-,\C)}(\C,U),\cl^2,\ccr^1)$ to be a right
$\C$-codual of $U$.\\

Consider a generic element $X\in^\da\! U$. Since $X:\C\ra U$ is a left
module map thus from (18)
\be
\cl\c X=(\id\ot X)\c\co .
\ee
Due to (13,14) the left and the right comultiplications on $X$ reads as follows
\be
\cl^1(X) = (X\ot\id)\c\co , \ \ \ \ \ \ccr^2(X) = \ccr\c X
\ee
modulo the identifications indicated in Remark 2.4.\medskip\\
{\bf Proposition 2.8}. \bm For a counital coassociative coalgebra $\C$
the bicomodules $^\da\!\C$ and $\C^\bu$ are isomorphic.\en

Instead of presenting the proof, we only notice that the required
isomorphism $^\da\!\C\cong \C^\bu$ can be realized by a mapping
\ $\tilde\ep :^\da\!\C\ni X\mapsto X^\ep\doteq\ep\c X\in\C^\bu$
induced by the counit $\ep$ (see Example 2.1 in this context).

\section{First order codifferential calculus}

According to the general scheme \cite{Abe}, the definition of a codifferential
calculus can be obtained by systematical reversing of all arrows  in the
diagrams defining the first order differential calculus on the algebra A.
One has also to replace the words: algebra, bimodule and multiplication
by its dual counterparts: coalgebra, bicomodule and comultiplication, etc.
An interesting notion of braided Hopf algebra {\it biderivation} has been
introduced in \cite{OPR}.\medskip\\
{\bf Definition 3.1}. A first order {\it codifferential calculus} on a
coalgebra $\C$ (FOCC in short)  consists of a $\C$-bicomodule $(U,\cl,\ccr)$
and a $\bk$-linear map (called a {\it coderivation}) $\de\in\hom(U,\C)$
such that
\be
\co\c\de = (\id\ot\de)\c\cl + (\de\ot\id)\c\ccr \ .
\ee

Writing down in the coordinate language $\de=e^k\ot\de_k$ with $\de_k\in\C$
uniquely defined for a given basis $\{e_k\}$ in $U$ one deduces from (21)
\be
\co (\de_k)= L_k^i\ot\de_i + \de_i\ot R^i_k
\ee
as the restriction on the coefficients $\de_k$. In particular, we can
consider a coderivation $\xi\in {\rm Coder}(\C)$ from the coalgebra
to itself. That is
\be
\co\c\xi = (\id\ot\xi + \xi\ot\id)\c\co \ .
\ee

Let us fix some FOCC $\de :U\ra\C$. In analogy with the algebra
case, the elements of $U$ can be called $1$-{\it form cofields}. Let us
take the left codual $^\da U$ of $U$: its elements can be called (left)
{\it vector cofields}. Now, with any vector cofield
$X\doteq X^i\ot e_i\in\,^\da U$
one can associate an endomorphism $X^\de\in\nd\C$
\be
X^\de\doteq \de\c X = X^i\ot\de_i .
\ee
This  becomes the dual counterpart of the famous Cartan formula
(cf. \cite{AB1,AB2}) Our task here is to investigate properties of
the map $\tilde\de :\,^\da\! U\ni X\mapsto X^\de\in\nd\C$.
This is done by the following\medskip\\
{\bf Theorem 3.2}. \bm
The map $\tilde\de :\,^\da\! U\ra\nd\C$ is a left comodule map.
Moreover, for any $X\in\,^\da\! U$ one has:
\be
\co\c X^\de = (\id\ot X^\de)\c\co\ +\ \ccr (X)^{(\de\ot\id)}\ .
\ee\en
{\it Proof:}\  The first claim requires the equality
$$
\cl^{{\small \nd}\C}\c\tilde\de = (\id\ot\tilde\de)\c\cl^{^\da\! U}
$$
which follows easily from (17) and (20).
To show the second, compose both sides of (21) with $X$. Thus
\be
\co\c X^\de=(\id\ot\de)\c\cl\c X + (\de\ot\id)\c\ccr\c X .\ee
From (19) one calculates $(\id\ot\de)\c\cl\c X=(\id\ot\de)\c(\id\ot X)
\c\co=(\id\ot X^\de)\c\co$. Applying (20) to the second term gives
$
(\de\ot\id)\c\ccr\c X= (\de\ot\id)\c\ccr^2(X)$.
This, obviously can be rewritten as $\ccr^2(X)^{(\de\ot\id)}$.
The proof is done.\hfill $\Box$\\

Because of the second term on the RHS in (25) and comparing with (23),
one can say that $X^\de$ is a deformed coderivation (see also \cite{AB1,AB2}
for a similar situation concerning vector fields ). We are going
to explain how to get an undeformed co-Leibniz rule (23).\\

Assume now that the coalgebra $\C$ is cocommutative. There exists
a well-known one-to-one correspondence between comodules and
cobimodules having the same left and right comultiplication, i.e.
$\ccr=\sw\c\cl$, where $\sw:\C\ot U\ra U\ot\C$ denotes the
canonical switch. In this case, one can recalculate the second
term in (26): $(\de\ot\id)\c\ccr\c X=(\de\ot\id)\c\sw\c\cl\c X =
(\de\ot\id)\c\sw\c(\id\ot X)\c\co =(X^\de\ot\id)\c\co$. This
proves \medskip\\ 
{\bf Theorem 3.3}. \bm Let $U$ be  a comodule and let $\de:U\ra\C$
be a FOCC over a cocommutative coalgebra $\C$.  Then for each
$X\in\,^\da\! U$ one has $X^\de\in {\rm Coder}(\C)$.\en\medskip\\

Let us finish by posing an unsolved question: Does there exist a canonical
(K\"ahler type) codifferential calculus $\de_0:U_0\ra\C$ on an
arbitrary cocommutative coalgebra
$\C$ such that $^\da\! U_0 = {\rm Coder}(\C)$ ?

\bigskip
{\small 
This work was supported by Polish KBN (grant \# 2 P03B 109 15) and
CONACyT, M\'exico (proyecto \# 27670 E). The authors would like to
thank Zbigniew Oziewicz and James Stasheff for their comments and
suggestions.}

\end{document}